\documentclass[smallcondensed]{svjour3}       % onecolumn (second format)
\smartqed  % flush right qed marks, e.g. at end of proof
\usepackage{graphicx}
\usepackage{hyperref}
\usepackage{url}
\usepackage{natbib}
\usepackage{subfigure}
\usepackage{tabularx}
\usepackage{mathtools}
\usepackage{amsmath}
\usepackage{amssymb}
\usepackage{algorithm}
\usepackage{algorithmic}
\mathtoolsset{showonlyrefs=true}
\newcommand{\argmax}{\operatorname*{argmax}} %\operatorname* pour les op. pouvant admettre des limites...
\newcommand{\argmin}{\operatorname*{argmin}}

\begin{document}

\title{Difference of Convex Functions Programming Applied to Control with Expert Data}
\titlerunning{DC Programming Applied to Control with Expert Data}   

\author{Bilal Piot$^{1}$ \\ Matthieu Geist$^{2}$ \\ Olivier Pietquin$^{3}$}

\institute{$^{1}$Univ. Lille, CRIStAL (UMR CNRS 9189/Lille 1) - SequeL team, bilal.piot@univ-lille1.fr \and $^{2}$CentraleSup\'elec (UMI 2958 GeorgiaTech-CNRS) - IMS-MaLIS team, matthieu.geist@centralesupelec.fr \and $^{3}$Univ. Lille, IUF, CRIStAL (UMR CNRS 9189/Lille 1) - SequeL team, olivier.pietquin@univ-lille1.fr, now with Google DeepMind}

\date{Received: date / Accepted: date}
\maketitle
\begin{abstract}
This paper reports applications of Difference of Convex functions (DC) programming to Learning from Demonstrations (LfD) and Reinforcement Learning (RL) with expert data. This is made possible because the norm of the Optimal Bellman Residual (OBR), which is at the heart of many RL and LfD algorithms, is DC. Improvement in  performance is demonstrated on two specific algorithms, namely Reward-regularized Classification for Apprenticeship Learning (RCAL) and Reinforcement Learning with Expert Demonstrations (RLED), through experiments on generic Markov Decision Processes (MDP), called Garnets.  
\end{abstract}

\section{Introduction}
\label{section: Introduction}
The Optimal Bellman Residual (OBR), being a core optimization criterion in Reinforcement Learning (RL), as been recently proven to be a Difference of Convex (DC) functions~\citep{piotNIPS2014}. As a consequence, this paper aims at extending previous results obtained in DC for batch RL to the fields of control with expert data and Learning from Demonstrations (LfD). More precisely, its objective is to leverage the knowledge in DC programming in order to improve the performance of existing methods in control and LfD.  

In control theory, there are two canonical ways to make an apprentice agent learn a task from expert demonstrations. The first one consists in directly learning a behaviour (mapping situations to decisions) generalising the expert decisions in states that were unvisited during demonstrations. This is the framework of LfD~\citep{pomerleau1989alvinn,atkeson1997robot,schaal1997learning,argall2009survey}. The second approach consists in inferring a goal that the apprentice agent should achieve from the demonstrations. The apprentice would than have to interact with the environment and find a strategy to attain it. When this goal is defined trough a reward function representing the local benefit of doing a particular action in a given state, the agent aims at maximising the sum of rewards encountered during its interaction with the environment. This is the framework of Reinforcement Learning (RL)~\citep{sutton1998reinforcement}. From a human perspective, the reward represents a local satisfaction and realising a task consists in maximising the sum of the local satisfactions. RL has the advantage of clearly defining a task through a reward function, providing with an optimisation criterion. However, optimising a sparsely distributed reward is sometimes a tricky task (sparse rewards are mainly encountered in practice because there are easier to define) as the agent has often to explore exhaustively its environment to discover rewards. For this reason, it can be useful to combine RL with LfD~\citep{schaal1997learning,Supelec885} in order to avoid learning from scratch and knowing precisely the task to achieve. There is a vast literature on LfD and RL but very few articles on how DC techniques can improve those methods (one exception being~\citep{piotNIPS2014}). 

DC programming~\citep{tao1997convex,tao2005dc} can transform a complex non-convex (but DC) problem into a series of simpler convex problems solvable via gradient-descent/ascent methods or Linear Programming (LP). This property is very interesting as it allows leveraging the huge amount of gradient-descent/ascent and LP literature. Thus, DC techniques have been applied to several domains and Machine Learning (ML) is no exception~\citep{le2014dc,le2014feature}. Indeed, DC methods can be used to address classification tasks~\citep{le2008dc,le2015sparse} and RL problems~\citep{piotNIPS2014}. 

%Therefore, we want to leverage the DC programming techniques developed in those previous papers and apply them to RL with expert data and LfD. 

To make use of DC, we place ourselves in the Markov Decision Process (MDP) paradigm, well-suited to study sequential decision making problems in stochastic, discrete-time and finite action-state space environments. In this specific framework, finding an optimal control can be cast into minimising a criterion, namely the OBR, which appears to be a DC function~\citep{piotNIPS2014}. More precisely, we focus on two existing methods called Reward-regularised Classification for Apprenticeship Learning (RCAL)~\citep{Supelec874} and Reinforcement Learning with Expert Demonstrations (RLED)~\citep{Supelec885}. RCAL is an LfD method which consists in constraining a classification method to using the dynamics of the underlying MDP so as to obtain better generalisation properties. RLED is an RL method using expert demonstrations together with OBR minimisation so as to boost learning speed. Those two algorithms consist in minimising a regularised classification criterion where the regularisation term, which is an empirical version of the OBR, is DC. In their original form, these algorithms where based on standard (sub-)gradient descent, but here we show how using DC techniques can improve the optimisation result.

The remaining of the paper is organised as follows. First, Sec.~\ref{section: Background} provides the notations and the background. It introduces the concepts of MDP, RL, IL, RLED and DC functions. Then, in Sec.~\ref{section: DC}, we show how RCAL and RLED can be decomposed in DC problems. Finally, in Sec~.\ref{section: Experiments}), we show experimental results.
\section{Background}
\label{section: Background}
Here, we introduce concepts such as MDP (Sec.~\ref{subsection: MDP}) and RL (Sec.~\ref{subsection: RL}) which are prerequisites to understanding LfD (Sec.\ref{subsection: LfD}) and RLED (Sec.\ref{subsection: RLED}) frameworks. We also briefly introduce the few basic notions of DC programming (Sec~\ref{subsection: DC}) required to decompose RCAL and RLED into DC problems.

Let us start with the general notations used throughout this paper. Let $(\mathbb{R},|.|)$ be the real space with its canonical norm and $X$ a finite set, $\mathbb{R}^X$ is the set of functions from $X$ to $\mathbb{R}$. The set of probability distributions over $X$ is noted $\Delta_X$. Let $Y$ be a finite set, $\Delta_X^Y$ is the set of functions from $Y$ to $\Delta_X$. Let $\alpha\in\mathbb{R}^X$, $p\geq1$ and $\nu\in\Delta_X$, we define the $\mathbf{L}_{p,\nu}$-semi-norm of $\alpha$, noted $\|\alpha\|_{p,\nu}$, by: $\|\alpha\|_{p,\nu}=(\sum_{x\in X}\nu(x)|\alpha(x)|^p)^{\frac{1}{p}}$. In addition, the infinite norm is noted $\|\alpha\|_{\infty}$ and defined as $\|\alpha\|_{\infty}=\max_{x\in X}|\alpha(x)|$. Let $v$ be a random variable which takes its values in $X$, $v\sim\nu$ means that the probability that $v=x$ is $\nu(x)$.
\subsection{Markov Decision Process}
\label{subsection: MDP}
The MDP paradigm is a state-of-the-art framework to learn optimal control in a stochastic, discrete-time and finite action-state space environment~\citep{howard1960dynamic,bellman1965dynamic}. Via the definition of a reward function representing the local information on the benefit of doing an action in a given state, it allows to find an optimal behaviour w.r.t. a predefined criterion. Here, we consider infinite-horizon $\gamma$-discounted MDPs where the chosen criterion is the $\gamma$-discounted and expected cumulative reward collected by an agent (a formal definition is given below). An optimal behaviour thus maximises the previous criterion and can be exactly found through Dynamic Programming (DP)~\citep{bertsekas1995dynamic,puterman1994markov} techniques such as Value Iteration (VI), Policy Iteration (PI) or Linear Programming (LP).

Here, the agent is supposed to act in a finite MDP represented by a tuple $M=\{S,A,R,P,\gamma\}$ where $S=\{s_i\}_{1\leq i\leq N_S}$ is the finite state space ($N_S\in\mathbb{N}^*$), $A=\{a_i\}_{1\leq i\leq N_A}$ is the finite action space ($N_A\in\mathbb{N}^*$), $R\in\mathbb{R}^{S\times A}$ is the reward function ($R(s,a)$ represents the local benefit of doing action $a$ in state $s$), $\gamma\in ]0,1[$ is a discount factor and $P\in \Delta_{S}^{S\times A}$is the Markovian dynamics which gives the probability, $P(s'|s,a)$, to reach $s'$ by choosing action $a$ in state $s$. It has been shown~\citep{puterman1994markov} for finite MDPs that it is sufficient to consider deterministic policies in order to obtain an optimal behaviour with respect to the $\gamma$-discounted and expected cumulative reward criterion. A deterministic policy $\pi$ is an element of $A^{S}$, it maps each state to a unique action and thus defines the behaviour of an agent. The quality of a policy $\pi$ is defined by the action-value function. For a given policy $\pi$, the action-value function $Q^\pi\in\mathbb{R}^{S\times A}$ is defined as:
\begin{equation}
Q^\pi(s,a)=\mathbb{E}^\pi[\sum_{t=0}^{+\infty}\gamma^tR(s_t,a_t)], 
\end{equation}
where $\mathbb{E}^\pi$ is the expectation over the distribution of the admissible trajectories $(s_0,a_0,s_1,\pi(s_1),\dots)$ obtained by executing the policy $\pi$ starting from $s_0=s$ and $a_0=a$. Therefore, the quantity $Q^\pi(s,a)$ represents the $\gamma$-discounted and expected cumulative reward collected by executing the policy $\pi$ starting from $s_0=s$ and $a_0=a$. Often, the concept of value function $V^\pi$ is used and corresponds to:
\begin{equation}
\forall s\in S,\quad V^\pi(s)=Q^\pi(s,\pi(s)),
\end{equation}
which represents the $\gamma$-discounted and expected cumulative reward collected by executing the policy $\pi$ starting from $s_0=s$ and $a_0=\pi(s)$. So, the aim, when optimising an MDP, is to find a policy $\pi$, called an optimal policy, such that:
\begin{equation}
\forall \pi'\in A^S,\forall s\in S,  V^\pi(s)\geq V^{\pi'}(s).
\end{equation}
To do so, an important tool is the optimal action-value function $Q^*\in\mathbb{R}^{S\times A}$ defined as $Q^*=\max_{\pi\in A^{S}}{Q^\pi}$. It has been shown by~\citet{puterman1994markov} that a policy $\pi$ is optimal if and only if $\forall s\in S,\quad V^\pi(s)=V^*(s)$, where $\forall s\in S,\quad V^*(s)=\max_{a \in A} Q^*(s,a)$. In addition, an important concept, that we will use throughout the paper, is greediness. A policy $\pi$ is said greedy with respect to a function $Q\in\mathbb{R}^{S\times A}$ if:
\begin{equation}
\forall s\in S, \pi(s)\in\argmax_{a\in A} Q(s,a).
\end{equation}
Greedy policies are important because a policy $\pi$ greedy with respect to $Q^*$ is optimal~\citep{puterman1994markov}. Thus, if we manage to find $Q^*$, we automatically found an optimal policy by taking a greedy policy with respect to $Q^*$. Moreover,
$Q^\pi$ and $Q^*$ are known to be the unique fixed points of the contracting operators $T^\pi$ and $T^*$ (also called Bellman operators) respectively:
\begin{align}
&\forall Q\in\mathbb{R}^{S\times A}, \forall (s,a)\in S\times A,\quad T^\pi Q(s,a)=R(s,a)+\gamma\sum_{s'\in S}P(s'|s,a)Q(s,\pi(s')),
\\
&\forall Q\in\mathbb{R}^{S\times A}, \forall (s,a)\in S\times A,\quad T^* Q(s,a)=R(s,a)+\gamma\sum_{s'\in S}P(s'|s,a)\max_{b\in A}Q(s,b).
\end{align}
This means, by uniqueness of the fixed points $Q^\pi$ and $Q^*$, that:
\begin{align}
Q^\pi=\argmin_{Q\in\mathbb{R}^{S\times A}}\|T^\pi Q-Q\|_{p,\mu}=\argmin_{Q\in\mathbb{R}^{S\times A}}\|T^\pi Q-Q\|_{\infty},
\\
Q^*=\argmin_{Q\in\mathbb{R}^{S\times A}}\|T^* Q-Q\|_{p,\mu}=\argmin_{Q\in\mathbb{R}^{S\times A}}\|T^* Q-Q\|_{\infty},
\label{eq:equation opti}
\end{align}
where $\mu\in\Delta_{S\times A}$ is such that $\forall (s,a)\in S\times A,\mu(s,a)>0$ and $p\geq1$. Thus, Eq~\eqref{eq:equation opti} shows that optimising an MDP can be seen as the minimisation of the criterion $J_{p,\mu}(Q)=\|T^* Q-Q\|_{p,\mu}$ where $T^*Q-Q$ is the OBR. Moreover, ~\citet{piotNIPS2014} showed that the function $J_{p,\mu}$ is DC and they provided an explicit decomposition for $p=1$ and $p=2$. However, minimising directly $J_{p,\mu}$ via a DC programming technique, such as DC Algorithm (DCA), when the MDP is perfectly known is not useful as there exists DP techniques (VI, PI and LP) which efficiently and exactly compute optimal policies. They rely on nice properties of the Bellman operators such as being a contraction and monotonicity. However, when the state space becomes large, two important problems arise and DP programming is not an option anymore. The first one, called the \textit{representation} problem, relate to the fact that  each value $Q^\pi(s,a)$ of the action-value functions cannot be stored as the number of values grows too much. So these functions need to be approximated with a moderate number of coefficients. The second problem, called the \textit{sampling} problem, arises because only samples from the Bellman operators are observed (there is thus only a partial knowledge of the dynamics $P$ and the reward function $R$). One solution for the representation problem is to use a linear approximation of the action-value functions thanks to a basis of $d\in\mathbb{N}^*$ functions $\phi=(\phi_i)_{i=1}^{d}$ where $\phi_i\in\mathbb{R}^{S\times A}$. In addition, we define for each state-action couple $(s,a)$ the vector $\phi(s,a)\in\mathbb{R}^d$ such that $\phi(s,a)=(\phi_i(s,a))_{i=1}^d$. Thus, the action-value functions are characterised by a vector $\theta\in\mathbb{R}^d$ and noted $Q_\theta$:
\begin{equation}
\forall \theta\in\mathbb{R}^d, \forall (s,a)\in S\times A, Q_{\theta}(s,a)=\sum_{i=1}^d\theta_i\phi_i(s,a)=\langle\theta,\phi(s,a)\rangle,
\end{equation}
where $\langle.,.\rangle$ is the canonical dot product of $\mathbb{R}^d$. In order to tackle the \textit{sampling} problem, RL techniques have been proposed and rely principally on Approximate DP (ADP) methods. One notable exception, developed by ~\citet{piotNIPS2014}, consists in directly minimising an empirical norm of the OBR via DCA (see Sec.~\ref{subsection: RL}).
\subsection{Reinforcement Learning}
\label{subsection: RL}

RL is a vast domain with different settings~\citep{sutton1998reinforcement,lange2012batch} sharing the property that the model (the dynamics $P$ and the reward $R$) of the MDP is only known through data (\textit{sampling} problem). Data can be collected on-line by an agent acting in the MDP, via a simulator or via previous interactions. Here, we focus on the last setting called the batch setting. More precisely, the machine is provided with a set of traces of interactions with the MDP:
\begin{equation}
D_{RL}=(s_j,a_j,r_j,s_j')_{j=1}^{N_{RL}}, 
\end{equation} 
where $s_j\in S$, $a_j\in A$, $r_j=R(s_j,a_j)$, $s_j'\sim P(.|s_j,a_j)$ and $N_{RL}\in\mathbb{N}^*$. Using only that set, a batch RL technique must estimate an optimal policy. Several techniques inspired by ADP such as Fitted-Q~\citep{ernst2005tree} (Approximate VI technique), Least squares Policy Iteration (LSPI)~\citep{lagoudakis2003least} (Approximate PI technique) and Locally Smoothed Regularised Approximate Linear Programming (LSRALP)-~\citep{taylor2012feature} exist. They rely on particular properties of the optimal Bellman operator $T^*$, such as monotonicity and contraction, to estimate the fixed-point $Q^*$. Here, we are interested in a new breed of RL techniques consisting in directly minimising the empirical OBR:
\begin{align}
J_{RL}(Q)&=\frac{1}{N_{RL}}\sum_{j=1}^{N_{RL}}|T^* Q(s_j,a_j)-Q(s_j,a_j) |,
\\
&=\frac{1}{N_{RL}}\sum_{j=1}^{N_{RL}}|R(s_j,a_j)+\gamma\sum_{s'\in S}P(s'|s_j,a_j)\max_{a\in A}Q(s',a)-Q(s_j,a_j) |.
\label{eq: RL def}
\end{align}
\citet{piotNIPS2014} showed that minimising directly the empirical optimal Bellman residual ($J_{RL}$) is a legit technique as:
\begin{equation}
\label{bound OBRM}
\|Q^*-Q^\pi\|_{p,\nu}\leq\frac{C}{1-\gamma}J_{RL}(Q),
\end{equation}
where $C$ is a constant depending on the dynamics $P$, $\nu\in\Delta_{S\times A}$ and $\pi$ is a greedy policy with respect to $Q$. This bound shows that minimising $J_{RL}$, leads to learn a function close to the optimal quality function $Q^*$.  A finite-sample analysis version of the bound and a comparison to ADP techniques is also provided by~\citet{piotNIPS2014}. In order to minimise $J_{RL}$, they showed that this criterion is DC, gave an explicit decomposition and proposed to use DCA. 

In practice, $T^*Q(s_j,a_j)$ cannot be computed from $D_{RL}$, as $P$ is unknown, unless the MDP is deterministic. Indeed, in that case, which is the one considered in our experiments (see Sec.~\ref{section: Experiments}), we have:
\begin{equation}
\sum_{s'\in S}P(s'|s_j,a_j)\max_{a\in A}Q(s',a)=\max_{a\in A}Q(s'_j,a).
\end{equation}
So, $T^*Q(s_j,a_j)=R(s_j,a_j)+\gamma\max_{a\in A}Q(s'_j,a)$ can be obtained from $D_{RL}$. In the general case, only a unbiased estimate of $T^*Q(S_i,A_i)$ can be computed via: 
\begin{equation}
\hat{T}^*Q(s_i,a_i)=R(s_i,a_i)+\gamma\max_{b\in A}Q(s'_i,b).
\end{equation}
The problem is that $|\hat{T}^*Q(s_j,a_j)-Q(s_j,a_j)|^p$ is a biased estimator of $|T^*Q(s_j,a_j)-Q(s_j,a_j)|^p$ and the bias is uncontrolled~\citep{antos2008learning}. In order to alleviate this typical problem, several better estimators $|\hat{T}^*Q(s_j,a_j)-Q(s_j,a_j)|^p$ of $|T^*Q(s_j,a_j)-Q(s_j,a_j)|^p$ have been proposed, such as embeddings in Reproducing Kernel Hilbert Spaces (RKHS)\citep{lever2012modelling} or locally weighted averager such as Nadaraya-Watson estimators\citep{taylor2012feature}. In both cases, the unbiased estimate of $T^*Q(s_j,a_j)$ takes the form:
\begin{equation}
\hat{T}^*Q(s_j,a_j)=R(s_j,a_j)+\gamma\frac{1}{N_{RL}}\sum_{j=k}^{N_{RL}}\beta_j(s'_k)\max_{a\in A}Q(s'_k,a), 
\end{equation}
where $\beta_j(s'_k)$ is the weight of samples $s'_k$ in the $T^*Q(s_j,a_j)$ estimate. 
%Thus, one can replace $\sum_{s'\in S}P(s'|s_j,a_j)\max_{a\in A}Q(s',a)$ by $\frac{1}{N}\sum_{k=1}^N\beta_j(s'_k)\max_{a\in A}Q(s'_k,a)$ to obtained an unbiased estimate of $J_{RL}$,.

\subsection{LfD and RCAL}
\label{subsection: LfD}
LfD consists in learning an policy $\pi_E$ from  demonstrations of an expert (which can be optimal or near-optimal) and without observing rewards. The aim is to generalise the expert behaviour in states that where not observed in the demonstration set. This setting is easily motivated as, in a lot of practical applications, it is easier to provide expert demonstrations than a reward function. Here, we consider the batch LfD setting where a set of expert demonstrations 
\begin{equation}
D_E=(s_i,a_i)_{i=1}^{N_E},
\end{equation} 
with  $s_i\in S$, $a_i=\pi_E(s_i)$ and $N_E\in\mathbb{N}^*$, and a set of sampled transitions without rewards
\begin{equation}
D_{NE}=(s_j,a_j,s_j')_{j=1}^{N_{NE}},
\end{equation}
with $s_j\in S$, $a_j\in A$, $s_j'\sim P(.|s_j,a_j)$ and $N_{NE}\in\mathbb{N}^*$, are provided. The set $D_{NE}$ gives a useful information on the dynamics of the underlying MDP and the set $D_E$ gives examples of an optimal (or sub-optimal) behaviour. There are several ways to tackle the LfD problem. The most well-known and studied are Inverse Reinforcement Learning (IRL)~\citep{ng2000algorithms,ziebart2008maximum,syed2008game,klein2012inverse,Supelec854} and Imitation Learning (IL)~\citep{ratliff2007imitation,syed2010reduction,ross2010efficient,ross2011reduction,judah2012active}. IRL consists in estimating a reward function that explains the expert behaviour. Once the reward is estimated, the resulting MDP has to be solved to end up with an actual policy. The interested reader can refer to this survey~\citep{neu2009training}.
On the other hand. IL consists in directly learning a mapping from states to actions to imitates the expert. This approach can be cast to a pure Supervised Learning (SL) problem such as Multi-Class Classification (MCC)~\citep{pomerleau1989alvinn,ratliff2007imitation,ross2010efficient,syed2010reduction}. Indeed, to compare to the standard classification notations, a state-action couple $(s_i,a_i=\pi_E(s_i))$ of the expert set $D_E$ could be seen as an input-label $(x_i,y_i)$ couple of a training set $D=(x_i\in X,y_i\in Y)_{i=1}^{N}$, where $X$ is a compact set of inputs (in our particular case $X$ is even finite) and $Y$ a finite set of labels. The goal of MCC is, given $D$, to find a decision rule $g\in\mathfrak{H}\subset Y^X$, where $\mathfrak{H}$ is an hypothesis space, that generalises the relation between inputs and labels by minimising  the empirical risk:
\begin{equation}
g=\argmin_{h\in\mathfrak{H}}\frac{1}{N}\sum_{i=1}^N \mathbf{1}_{\{y_i=h(x_i)\}}.
\end{equation}
Properties of $g$ and how well it generalises the data are notably studied by~\citet{vapnik1998statistical}. However, minimising directly the empirical risk is unrealistic and practitioners use convex surrogates. Often, another approach, called score-based MCC, where a score function $Q\in\mathbb{R}^{X\times Y}$ is learnt, is used (we intentionally decide to use the same notation as action-value functions as there is a close link between between score functions and action-value functions~\citep{Supelec854}). The score $Q(x,y)$ represents the correspondence between the input $x$ and the label $y$. The higher the score is, the more likely $y$ will be chosen when $x$ is the input. Thus, the decision rule $g$ corresponding to the score $Q$ is $g(x)=\argmax_{y\in Y}Q(x,y)$. For instance,~\citet{ratliff2007imitation} use a large margin approach which is a score-based MCC for solving an IL problem. The large margin approach consists, given the training set $D$, in minimising the following criterion:
\begin{align}
\label{equation : structured margin function}
J(Q)&=\frac{1}{N}\sum_{i=1}^{N}\max_{y\in Y}\{Q(x_i,y)+l(x_i,y_i,y)\}-Q(x_i,y_i),
\end{align}
where $l\in\mathbb{R}_+^{X\times Y\times Y}$ is called the margin function. If this function is zero, minimising $J(Q)$ attempts to find a score function for which the example labels are scored higher than all other labels. Choosing a non-zero margin function improves generalisation~\citep{ratliff2007imitation} and instead of requiring only that the example label is scored higher than all other labels, it requires it to be better than each label $y$ by an amount given by the margin function. In practice, one can use a margin function equals $0$ for the inputs $(x_i,y_i,y=y_i)$ and $1$ otherwise (this is the margin function we chose in our experiments). Applying the large margin approach to the LfD problem gives the following minimisation criterion:
\begin{align}
\label{equation : structure margin MDP}
J_E(Q)&=\frac{1}{N_E}\sum_{i=1}^{N_E}\max_{a\in A}\{Q(s_i,a)+l(s_i,\pi_E(s_i),a)\}-Q(s_i,\pi_E(s_i)).
\end{align}
However this approach does not take into account the underlying dynamics of the MDP represented in the set $D_{NE}$.

To avoid that drawback,~\citet{Supelec874} propose to see the score function $Q$ as an optimal quality function $Q^*$ of an MDP. To do so, they rely on the one-to-one relation between optimal quality functions and rewards functions. Indeed, for each function $Q\in\mathbb{R}^{S\times A}$, there exists a reward function $R_Q\in\mathbb{R}^{S\times A}$ such that $Q=Q^*$ where $Q^*$ is the optimal quality function with respect to the reward $R_Q$. Moreover, there is an explicit formula for $R_Q$ depending only on $Q$ and $P$~\citep{Supelec874}:
\begin{equation}
R_Q(s,a)=Q(s,a)-\gamma \sum_{s'\in S}P(s'|s,a)\max_{b\in A}Q(s',b).
\end{equation}
Knowing that,~\citet{Supelec874} propose to regularise the criterion $J_E$ by a term controlling the sparsity of the reward associated to the score function. This regularisation term is:
\begin{align}
J_{NE}(Q)&=\frac{1}{N_{NE}}\sum_{j=1}^{N_{NE}}|R_Q(s_j,a_j)|,
\\
&=\frac{1}{N_{NE}}\sum_{j=1}^{N_{NE}}|\gamma\sum_{s'\in S}P(s'|s_j,a_j)\max_{a\in A}Q(s',a)-Q(s_j,a_j) |.
\label{eq: NE def}
\end{align}
This helps to reduce the variance of the method as it considers as good candidates only $Q$ functions with sparse rewards $R_Q$. The algorithm RCAL consists in minimising by a gradient descent the following criterion: 
\begin{equation}
J_{RCAL}(Q)=J_E(Q)+\lambda_{RCAL}J_{NE}(Q).
\end{equation}
However, it is easy to see that $J_{NE}(Q)$ (see Eq.~\eqref{eq: NE def}) corresponds to $J_{RL}(Q)$ (see Eq.~\eqref{eq: RL def}) when the reward function is null. Thus, $J_{NE}(Q)$ is also DC and as $J_E$ is convex, then $J_{RCAL}$ is DC. So, we propose to use the DCA to minimise $J_{RCAL}$ in Sec.~\ref{section: DC}.

\subsection{Reinforcement Learning with Expert Demonstrations}
\label{subsection: RLED}
RLED aims at finding an optimal control in a MDP where some expert data are provided in addition to standard sampled transitions with rewards. Such a paradigm is also easily motivated as in a lot of practical applications a goal (reward function) is provided to an agent (a robot for instance) but it can be quite difficult or risky to optimise from scratch (huge or dangerous environment to explore). Also a good control is often difficult to find as the reward function is very sparse and the agent needs to explore a lot of possibilities before finding a reward and retro-propagate it. Thus, an expert (a human for instance) can provide some demonstrations in order to guide the agent through the good learning path and accelerate the learning process~\citep{clouse1996introspection,gil2009cognitive,knox2012reinforcement,taylor2011integrating,griffith2013policy}. This combination of reward and expert data is somehow what we can experience in our daily life when we set goals to achieve (reward functions) and we observe other human beings (experts) achieving that same goals. Here, we consider the batch setting~\citep{kim2013learning,Supelec885}. More precisely, the apprentice agent is given a set of expert demonstrations (the same as the one in LfD)
\begin{equation}
D_E=(s_i,a_i)_{i=1}^{N_E},
\end{equation}
where $s_i\in S$, $a_i=\pi_E(s_i)$ and $N_E\in\mathbb{N}^*$, and a set of sampled transitions (the same as the one in RL)
\begin{equation}
D_{RL}=(s_j,a_j,r_j,s_j')_{j=1}^{N_{RL}},
\end{equation}
where $s_j\in S$, $a_j\in A$, $r_j=R(s_j,a_j)$, $s_j'\sim P(.|s_j,a_j)$ and $N_{RL}\in\mathbb{N}^*$. \citet{Supelec885} propose the RLED algorithm, minimising the following criterion:
\begin{equation}
J_{RLED}(Q)= J_E(Q) +\lambda_{RLED} J_{RL} (Q),
\end{equation}
combining two criteria: $J_E$ and $J_{RL}$ (defined in Eq.~\eqref{eq: RL def}). The regularisation factor $\lambda_{RLED}$ weights the importance between the expert and RL data. If one has a high confidence on the quality of the RL data, one will set  $\lambda_{RLED}$ to high value and to a low value otherwise. The criterion $J_{RLED}$ can also be seen as the minimisation of $J_{RL}$ guided by constraints provided by the expert data~\citep{Supelec885}. Another explanation could be that RLED produces a score function $Q$ that is forced to be an action-value function. This accelerates the optimisation of $J_{RL}$ and improves the performance of the method. In the original paper~\citep{Supelec885}, the authors propose to minimise $J_{RLED} (Q)$ by a gradient descent. However, as $J_{RL}$ is DC and $J_E$ is convex, then $J_{RLED}$ is DC. Thus, we propose to use DCA to minimise this criterion in Sec.~\ref{section: DC}. But before, we give some basics on DC programming which are sufficient to derive DC decompositions for RLED and RCAL.

\subsection{Basics on DC and DC programming}
\label{subsection: DC}
DC programming addresses non-convex (but DC) and non-differentiable (but sub-differentiable) optimisation problems by transforming them into a series of intermediary convex (thus simpler) optimisations problems. It allows leveraging the knowledge on convex optimisation and for that reason as been adapted to different domains such as Machine Learning~\citep{le2008dc,le2014dc,le2014feature}. It also gives some guarantees when one of the function of the DC decomposition is polyhedral, such as convergence in finite time to local minima. Thus, it seems a better solution than a simple gradient descent when confronted to complex non-convex (but DC) optimisation problems.

Let $E$ be a finite dimensional Hilbert space, $\langle .,.\rangle_E$ and $\|.\|_E$ its dot product and norm respectively. We say that a function $J\in \mathbb{R}^E$ is DC if there exists $f,g\in\mathbb{R}^E$ which are convex and lower semi-continuous such that $J=f-g$~\citep{tao2005dc}. The set of DC functions is noted $DC(E)$ and is stable to most of the operations that can be encountered in optimisation, contrary to the set of convex functions. Indeed, let $(J_i)_{i=1}^{K}$ be a sequence of $K\in\mathbb{N}^*$ DC functions and $(\alpha_i)_{i=1}^K\in\mathbb{R}^K$ then $\sum_{i=1}^K\alpha_iJ_i$, $\prod_{i=1}^KJ_i$, $\min_{1\leq i\leq K}J_i$, $\max_{1\leq i\leq K}J_i$ and $|J_i|$ are DC functions~\citep{hiriart1985generalized}. In order to minimise a DC function $J=f-g$, we need to define a notion of differentiability for convex and lower semi-continuous functions. Let $g$ be such a function and $e\in E$, we define the sub-gradient $\partial_e g$ of $g$ in $e$ as:
\begin{equation}
\partial_e g=\{\delta_e\in E, \forall e'\in E, g(e')\geq g(e)+\langle e'-e,\delta_e\rangle_E\}.
\end{equation}
For convenience, we make this little abuse of notations where $\partial_e g$ can refer to any element of $\partial_e g$. For a convex and lower semi-continuous $g\in \mathbb{R}^E$, the sub-gradient $\partial_e g$ is non empty for all $e\in E$~\citep{hiriart1985generalized}. This observation leads to a minimisation method of a function $J\in DC(E)$ called Difference of Convex functions Algorithm (DCA). Indeed, as $J$ is DC, we have:
\begin{equation}
\forall (e,e')\in E^2, J(e')=f(e')-g(e')\underset{(a)}{\leq} f(e')-g(e)-\langle e'-e,\partial_e g\rangle_E,
\end{equation}
where  inequality $(a)$ is true by definition of the sub-gradient. Thus, for all $e\in E$, the function $J$ is upper bounded by a function $I_e\in\mathbb{R}^E$ defined, $\forall e'\in E$, by 
\begin{equation}
I_e(e')=f(e')-g(e)-\langle e'-e,\partial_e g\rangle_E.
\end{equation}
The function $I_e$ is a convex and lower semi-continuous function (as it is the sum of two convex and lower semi-continuous functions which are $f$ and the linear function $\forall e'\in E,\langle e-e',\partial_e g\rangle_E-g(e)$). In addition, those functions have the particular property that $\forall e\in E, J(e)=I_e(e)$. The set of convex functions $(I_e)_{e\in E}$ that upper-bound the function $J$ plays a key role in DCA.

The algorithm DCA~\citep{tao2005dc} consists in constructing a sequence $(e_n)_{n\in\mathbb{N}}$ such that the sequence $(J(e_n))_{n\in\mathbb{N}}$ decreases. The first step is to choose a starting point $e_0\in E$, then to minimise the convex function $I_{e_0}$ that upper-bounds the function $J$. We can remark that minimising $I_{e}$ is equivalent to minimising $I'_e$ defined by  $\forall e'\in E$ 
\begin{equation}
I'_e(e')=f(e')-\langle e',\partial_e g\rangle_E.
\end{equation}
We note $e_1$ a minimiser of $I_{e_0}$, $e_1\in\argmin_{e\in E} I_{e_0}$. This minimisation can be realised by any convex optimisation solver. As $J(e_0)=I_{e_0}(e_0)\geq I_{e_0}(e_1)$ and $I_{e_0}(e_1)\geq J(e_1)$, then $J(e_0)\geq J(e_1)$. Thus, if we construct the sequence $(e_k)_{k\in\mathbb{N}}$ such that $\forall k\in\mathbb{N}, e_{k+1}\in\argmin_{e\in E} I_{e_k}$ and $e_0\in E$, then we obtain a decreasing sequence $(J(e_k))_{k\in\mathbb{N}}$. Therefore, the algorithm DCA solves a sequence of convex optimisation problems in order to solve a DC optimisation problem. Three important choices can radically change the DCA performance: the first one is the explicit choice of the decomposition of $J$, the second one is the choice of the starting point $e_0$ and finally the choice of the intermediate convex solver. The DCA algorithm hardly guarantees convergence to the global optima, but it usually provides good solutions. Moreover, it has some nice properties when one of the functions $f$ or $g$ is polyhedral. A function $g$ is said polyhedral when $\forall e\in E, g(e)=\max_{1\leq i \leq K} [\langle \alpha_i,e \rangle_H+\beta_i]$, where $(\alpha_i)_{i=1}^K\in E^K$ and $(\beta_i)_{i=1}^K\in\mathbb{R}^K$. If one of the function $f,g$ is polyhedral, $J$ is under bounded, the DCA sequence $(e_k)_{k\in\mathbb{N}}$ is bounded and the DCA algorithm converges in finite time to one of the local minima. The finite time aspect is important in terms of implementation. More details about DC programming and DCA are given by~\citet{tao2005dc} and even conditions for convergence to the global optima.

To summarise, once a DC decomposition is found, minimising the DC criterion $J=f-g$ via DCA corresponds to minimise the following intermediary convex functions:
\begin{equation}
I'_{k}(e')=f(e')-\langle e',\partial_{e_k}g\rangle_E,
\end{equation}
with $e_0\in E$ and $e_{k+1}\in\argmin_{e'\in E} I'_{k}(e')$. In practice, DCA stops when $e_{k+1}=e_{k}$ or when $k=K$ (with $K$ the maximal number of steps for DCA and this is the stopping criterion chosen in our experiments) and the output of the algorithm is $e_{k}$. In addition, obtaining $e_{k+1}$ can be done by linear programming (if $I'_{k}$ can be transformed into a linear program) or by gradient descent. In our experiments, we choose gradient descent to solve the intermediary convex problems with the following updates:
\begin{align}
e'_{0}=e_k,\quad \partial_{e'_{p}} I'_{k}=\partial_{e'_{p}}f-\partial_{e_k}g,\quad e'_{p+1}=e'_{p}-\alpha_{p}\frac{\partial_{e'_{p}} I'_{k} }{\|\partial_{e'_{p}} I'_{k}\|_E},
\label{eq: updategrad}
\end{align}
where $(\alpha_{p})_{p\in\mathbb{N}}\in\mathbb{R}_+^{\mathbb{N}}$. Finally, we set $e_{k+1}=e'_{p^*}$ where $p^*$ meets a stopping criterion such as $p^*=N$ (with $N$ is the maximal number of steps of the gradient descent and this is the stopping criterion chosen in our experiments) or $\partial_{e'_{p}} I'_{k}=0$ for instance.
Thus, when the intermediary convex problems $I'_{k}$ are determined, it is necessary to be able to compute their gradient $\partial_{e'_{p}} I'_{k}$ in order to apply DCA. In the next section, we give the decompositions of the criteria $J_{RCAL}$ and $J_{RLED}$ and how to compute the different gradients.
\section{DC Decompositions for RCAL and RLED }
\label{section: DC}
In this section, we derive a DC decomposition for the criteria $J_{RLED}$ (Sec.~\ref{subsection: RLED DC}) and $J_{RCAL}$ (Sec.~\ref{subsection: RCAL DC}) from the DC decompositions of $J_{RL}$ and $J_{NE}$ (Sec.~\ref{subsection: RL DC}). Several decompositions are possible, we describe the one that we actually use in experiments. Here, the DC decompositions is realised as if we could compute $\gamma\sum_{s'\in S}P(s'|s_j,a_j)\max_{a\in A}Q(s',a)$. In practice, in the deterministic case, we replace this quantity by $\gamma P(s_j'|s_j,a_j)\max_{a\in A}Q(s_j',a)$ which is easily computable and in the general case by $\frac{1}{N}\sum_{k=1}^N\beta_j(s'_k)\max_{a\in A}Q(s'_k,a)$ which is obtained using RKHS embedding or Nadaraya-Watson estimators.
\subsection{DC decomposition of $J_{RL}$ and $J_{NE}$}
\label{subsection: RL DC} 
Let us start with the criterion $J_{RL}(Q)$
\begin{align}
J_{RL}(Q)&=\frac{1}{N_{RL}}\sum_{j=1}^{N_{RL}}|T^* Q(s_j,a_j)-Q(s_j,a_j) |,
\\
&=\frac{1}{N_{RL}}\sum_{j=1}^{N_{RL}}|R(s_j,a_j)+\gamma\sum_{s'\in S}P(s'|s_j,a_j)\max_{a\in A}Q(s',a)-Q(s_j,a_j) |.
\end{align}
As we have seen previously, in RL, the functions $Q$  are characterised by a vector $\theta\in\mathbb{R}^d$ and noted $Q_{\theta}(s,a)=\sum_{i=1}^d\theta_i\phi_i(s,a)=\langle\theta,\phi(s,a)\rangle$. Thus, we consider the criterion:
\begin{equation}
J_{RL}(\theta)=\frac{1}{N_{RL}}\sum_{j=1}^{N_{RL}}|R(s_j,a_j)+\gamma\sum_{s'\in S}P(s'|s_j,a_j)\max_{a\in A}\langle \phi(s',a) , \theta \rangle -\langle \phi(s_j,a_j), \theta \rangle |.
\end{equation}
Noticing that $\gamma\sum_{s'\in S}P(s'|s_j,a_j)\max_{a\in A}\langle \phi(s',a) , \theta \rangle$ is convex in $\theta$ as a sum of a max of convex functions, that $\langle \phi(s_j,a_j), \theta \rangle$ is also convex in $\theta$ and that $|f-g|=2\max(f,g)-(f+g)$, we have the following DC decomposition for $J_{RL}$ (a complete proof is given by~\citet{piotNIPS2014}):
\begin{align}
f_{RL}^j(\theta)&= 2\max\left( R(s_j,a_j)+\gamma\sum_{s'\in S}P(s'|s_j,a_j)\max_{a\in A}\langle \phi(s',a) , \theta \rangle  , \langle \phi(s_j,a_j), \theta \rangle  \right),
\\
g_{RL}^j(\theta)&= R(s_j,a_j)+\gamma\sum_{s'\in S}P(s'|s_j,a_j)\max_{a\in A}\langle \phi(s',a) , \theta \rangle   + \langle \phi(s_j,a_j), \theta \rangle ,
\\
f_{RL}(\theta)&=\frac{1}{N_{RL}}\sum_{j=1}^{N_{RL}}f_{RL}^j(\theta),\quad g_{RL}(\theta)=\frac{1}{N_{RL}}\sum_{j=1}^{N_{RL}}g_{RL}^j(\theta),
\\
J_{RL}(\theta)&=\frac{1}{N_{RL}}\sum_{j=1}^{N_{RL}}f_{RL}^j(\theta) - g_{RL}^j{\theta}=f_{RL}(\theta)-g_{RL}(\theta).
\end{align}
We can do exactly the same calculus for $J_{NE}$ as it is the same criterion than $J_{RL}$ with a null reward. We have:
\begin{align}
J_{NE}(Q)&=\frac{1}{N_{NE}}\sum_{j=1}^{N_{NE}}|R_Q(s_j,a_j) |,
\\
J_{NE}(Q)&=\frac{1}{N_{NE}}\sum_{j=1}^{N_{NE}}|\gamma\sum_{s'\in S}P(s'|s_j,a_j)\max_{a\in A}Q(s',a)-Q(s_j,a_j) |.
\end{align}
With the linear parametrisation, we obtain: 
\begin{align}
J_{NE}(\theta)&=\frac{1}{N_{RL}}\sum_{i=j}^{N_{RL}}|\gamma\sum_{s'\in S}P(s'|s_j,a_j)\max_{a\in A}\langle \phi(s',a) , \theta \rangle -\langle \phi(s_j,a_j), \theta \rangle|,
\\
f_{NE}^j(\theta)&= 2\max\left(\gamma\sum_{s'\in S}P(s'|s_j,a_j)\max_{a\in A}\langle \phi(s',a) , \theta \rangle   , \langle \phi(s_j,a_j), \theta \rangle  \right),
\\
g_{NE}^j(\theta)&=\gamma\sum_{s'\in S}P(s'|s_j,a_j)\max_{a\in A}\langle \phi(s',a) , \theta \rangle   + \langle \phi(s_j,a_j), \theta \rangle,
\\
f_{NE}(\theta)&=\frac{1}{N_{NE}}\sum_{j=1}^{N_{NE}}f_{NE}^j(\theta),\quad g_{NE}(\theta)=\frac{1}{N_{NE}}\sum_{j=1}^{N_{NE}}g_{NE}^j(\theta),
\\
J_{NE}(\theta)&=\frac{1}{N_{NE}}\sum_{j=1}^{N_{NE}}f_{NE}^j(\theta) - g_{NE}^j(\theta)=f_{NE}(\theta)-g_{NE}(\theta).
\end{align}
Now that we have the DC decompositions, it is sufficient to calculate the intermediary convex problems ($I^k_{RL}(\theta)$ and $I^k_{NE}(\theta)$). To do so, we need the gradients $\partial_{\theta '} g_{RL}$ and $\partial_{\theta '} g_{NE}$:
\begin{align}
\partial_{\theta '} g_{RL}^j = \gamma\sum_{s'\in S}P(s'|s_j,a_j)\phi(s',a^*_{\theta',s'})   +\phi(s_j,a_j),\quad \partial_{\theta '} g_{RL}&=\frac{1}{N_{RL}}\sum_{j=1}^{N_{RL}}\partial_{\theta'} g_{RL}^j,
\end{align}
\begin{align}
\partial_{\theta '} g_{NE}^j =\gamma\sum_{s'\in S}P(s'|s_j,a_j)\phi(s',a^*_{\theta',s'}) +\phi(s_j,a_j),\quad \partial_{\theta '} g_{NE} =\frac{1}{N_{NE}}\sum_{j=1}^{N_{NE}}\partial_{\theta '} g_{NE}^j.
\end{align}
where $a^*_{\theta',s'}=\argmax_{a\in A} \langle \phi(s',a), \theta ' \rangle$. So the intermediary convex problems are:
\begin{align}
I^k_{RL}(\theta)=f_{RL}(\theta)-\langle \partial_{\theta^k_{RL}} g_{RL},\theta \rangle,\quad I^k_{NE}(\theta)=f_{NE}(\theta)-\langle\partial_{\theta^k_{NE}} g_{NE},\theta \rangle,
\end{align}
with $\theta^0_{RL}$, $\theta^0_{RL}$ in $\mathbb{R}^d$, $\theta^{k+1}_{RL}=\argmin_{\theta\in\mathbb{R}^d}I^k_{RL}(\theta)$ and $\theta^{k+1}_{NE}=\argmin_{\theta\in\mathbb{R}^d}I^k_{NE}(\theta)$.
To minimise $I^k_{RL}(\theta)$ and $I^k_{NE}(\theta)$, we have to compute their gradients and do the update as in Eq~.\eqref{eq: updategrad}:
\begin{align}
&\partial_{\theta'}f_{RL}^j=\left\{\begin{array}{l}
 \gamma\sum_{s'\in S}P(s'|s_j,a_j)\phi(s',a^*_{\theta',s'}) \\
 \text{ if $R(s_j,a_j)+\gamma\sum_{s'\in S}P(s'|s_j,a_j)\max_{a\in A}\langle \phi(s',a) , \theta' \rangle> \langle \phi(s_j,a_j), \theta' \rangle$},  \\
\phi(s_j,a_j) \text{ else }  
\end{array}\right.
\\
&\partial_{\theta'}f_{RL}=\frac{1}{N_{RL}}\sum_{j=1}^{N_{RL}}\partial_{\theta'}f_{RL}^j,\quad \partial_{\theta'} I^k_{RL}=\partial_{\theta'}f_{RL}(\theta)-\partial_{\theta^k_{RL}} g_{RL}.
\\
&\partial_{\theta'}f_{NE}^j=\left\{\begin{array}{l}
 \gamma\sum_{s'\in S}P(s'|s_j,a_j)\phi(s',a^*_{\theta',s'}) \\
 \text{ if $\gamma\sum_{s'\in S}P(s'|s_j,a_j)\max_{a\in A}\langle \phi(s',a) , \theta' \rangle> \langle \phi(s_j,a_j), \theta' \rangle$},  \\
\phi(s_j,a_j) \text{ else },  
\end{array}\right.
\\
&\partial_{\theta'}f_{NE}=\frac{1}{N_{NE}}\sum_{j=1}^{N_{NE}}\partial_{\theta'}f_{NE}^j,\quad
\partial_{\theta'} I^k_{NE}=\partial_{\theta'}f_{NE}(\theta)-\partial_{\theta^k_{RL}} g_{NE}.
\end{align}
The DC decompositions of $J_{RCAL}$ and $J_{RLED}$ follows directly from the ones of $J_{RL}$ and $J_{NE}$.
In addition, one can easily notice that $f_{RL}$, $g_{RL}$, $f_{NE}$ and $g_{NE}$ are polyhedral and that property will be directly transmitted to the decompositions of $J_{RCAL}$ and $J_{RLED}$.
\subsection{DC decomposition of $J_{RCAL}$}
\label{subsection: RCAL DC}
The criterion $J_{RCAL}$ is composed of two criterion $J_E$ and $J_{NE}$:
\begin{align}
J_E(Q)&=\frac{1}{N_E}\sum_{i=1}^{N_E}\max_{a\in A}[Q(s_i,a)+l(s_i,a_i,a)]-Q(s_i,a_i),
\\
J_{RCAL}(Q)&=J_E(Q)+\lambda_{RCAL}J_{NE}(Q).
\end{align}
With a linear parametrisation, we have:
\begin{align}
J_E(\theta)&=\frac{1}{N_E}\sum_{i=1}^{N_E}\max_{a\in A}[\langle \phi(s_i,a), \theta \rangle + l(s_i,a_i,a)] -\langle \phi(s_i,a_i), \theta\rangle,
\\
J_{RCAL}(\theta)&=J_E(\theta)+\lambda_{RCAL}J_{NE}(\theta).
\end{align}
Thus, the DC decomposition is quite trivial to obtain as $J_E$ is convex:
\begin{align}
f_{RCAL}(\theta)&=J_E(\theta)+\lambda_{RCAL}f_{NE}(\theta),
\\
g_{RCAL}(\theta)&=\lambda_{RCAL}g_{NE}(\theta),
\\
J_{RCAL}(\theta)&=f_{RCAL}(\theta)-g_{RCAL}(\theta).
\end{align}
To obtain the intermediary convex problems, we need to calculate the gradient of $g_{RCAL}$:
\begin{align}
\partial_{\theta '} g_{RCAL}&=\lambda_{RCAL} \partial_{\theta '}g_{NE},
\\
&=\frac{\lambda_{RCAL}}{N_{NE}}\sum_{j=1}^{N_{NE}}\partial_{\theta '} g_{NE}^j.
\end{align}
Thus, the convex intermediary problems have the following form:
\begin{align}
I^k_{RCAL}(\theta)&=f_{RCAL}(\theta)-\langle \partial_{\theta^k_{RCAL}} g_{RCAL},\theta\rangle.
\end{align}
To minimise $I^k_{RCAL}$, we calculate its gradient:
\begin{align}
\partial_{\theta'} J_E&=\frac{1}{N_E}\sum_{i=1}^{N_E}\phi(s_i,a_{i,\theta'}^*)-\phi(s_i,a_i),
\\
\partial_{\theta'} f_{RCAL}&=\partial_{\theta'} J_E + \lambda_{RCAL}\partial_{\theta'} f_{NE},
\\
\partial_{\theta'} I^k_{RCAL}&=\partial_{\theta'} f_{RCAL} -\partial_{\theta^k_{RCAL}} g_{RCAL}.
\end{align}
where $a_{i,\theta'}^*=\argmax_{a\in A}[\langle \phi(s_i,a), \theta' \rangle + l(s_i,a_i,a)]$.
\subsection{DC decomposition of $J_{RLED}$}
\label{subsection: RLED DC} 
The decomposition of $J_{RLED}$ is quite similar as the one of $J_{RCAL}$. We present it briefly for sake of completeness. $J_{RLED}$ is composed by two terms:
\begin{equation}
J_{RLED}(Q)=J_E(Q)+\lambda_{RLED}J_{RL}(Q).
\end{equation}
With a linear parametrisation, we have:
\begin{equation}
J_{RLED}(\theta)=J_E(\theta)+\lambda_{RLED}J_{RL}(\theta).
\end{equation}
As $J_E$ is convex, a DC decomposition is quite trivial to obtain: 
\begin{align}
f_{RLED}(\theta)&=J_E(\theta)+\lambda_{RLED}f_{RL}(\theta),
\\
g_{RLED}(\theta)&=\lambda_{RLED}g_{RL}(\theta),
\\
J_{RLED}(\theta)&=f_{RLED}(\theta)-g_{RLED}(\theta).
\end{align}
Now, to minimise $J_{RLED}$, we calculate the intermediary convex problems by obtaining the gradient of $g_{RLED}$:
\begin{align}
\partial_{\theta '} g_{RLED}&=\lambda_{RLED} \partial_{\theta '}g_{RL},
\\
&=\frac{\lambda_{RLED}}{N_{RL}}\sum_{j=1}^{N_{RL}}\partial_{\theta '} g_{RL}^j.
\end{align}
Thus, the intermediary convex problems have the following form:
\begin{align}
I^k_{RLED}(\theta)&=f_{RLED}(\theta)-\langle\partial_{\theta^k_{RLED}} g_{RLED},\theta\rangle.
\end{align}
Finally, in order to minimise $I^k_{RLED}$ by gradient descent, we need to calculate $\partial_{\theta'} I^k_{RLED}$:
\begin{align}
\partial_{\theta'} f_{RLED}&=\partial_{\theta'} J_E + \lambda_{RLED}\partial_{\theta'} f_{RL},
\\
\partial_{\theta'} I^k_{RLED}&=\partial_{\theta'} f_{RLED} -\partial_{\theta^k_{RLED}} g_{RLED}.
\end{align}
Now that we have the DC decompositions of $J_{RCAL}$ and $J_{RLED}$, we can compare the performance of those algorithms when the minimisation is realised via direct gradient descent or via DCA. This comparison is realised, in Sec.~\ref{section: Experiments}, on an abstract but representative class of MDPs called Garnets.
\section{Experiments}
\label{section: Experiments}
This section is composed of three experiments which aim at showing that using DC programming instead of gradient descent slightly improve the performance of existing algorithms, namely RCAL and RLED. The first experiment consists in showing the performance improvement of the RCAL algorithm when the set $D_{NE}$ is fixed and the set $D_{E}$ is growing on different MDPs which are randomly generated and called Garnets. The RCAL algorithm which consists in the minimisation of the criterion $J_{RCAL}$ is done by two methods. The first one is by gradient descent and is called RCAL in the remaining. The second one is by DCA and is called RCALDC. We also compare RCAL to a classical classification algorithm~\citep{ratliff2007imitation}, called Classif and which corresponds to the minimisation of $J_{RCAL}$ by gradient descent when $\lambda_{RCAL}=0$. The second experiments focuses on the performance improvement of RLED when the set $D_{RL}$ is fixed  and the set $D_{E}$ is growing. We compare the algorithm RLED when the minimisation is done by gradient descent, called RLED, and when the minimisation is done by DCA, called RLEDDC. Those algorithms are compared to LSPI which uses only the set $D_{RL}$ as input and Classif which uses only the set $D_E$ as input. Finally, the last experiment focuses on the performance improvement of RLED when the set $D_E$ is fixed and the set $D_{RL}$ is growing. But first, to realise those experiments, we need to introduce the notion of Garnets, explain how we construct the sets $D_E$, $D_{NE}$ and $D_{RL}$ and give the value of the different parameters of DCA and the gradient descent algorithms.

Garnets~\citep{archibald1995generation} are an abstract class of finite MDPs and easy to build. Here, we consider a special case of Garnets specified by three parameters: $(N_S,N_A,N_B)$. Parameters $N_S$ and  $N_A$ are respectively the number of states and actions. Thus, $S=(s_i)_{i=1}^{N_S}$ and $A=(a_i)_{i=1}^{N_A}$ are, respectively, the state and action spaces. The parameter $N_B$ ($N_B\leq N_S$), called the branching factor, defines for each state-action couple $(s,a)$ the number of next states. Here, we consider deterministic MDPs with $N_B=1$. The next state for each $(s,a)$, noted $s'_{s,a}$, is drawn uniformly from the set of states. In addition the discount factor is set to $\gamma=0.9$ or $\gamma=0.99$. Finally, we need to define the reward function $R$. To do so, we draw uniformly and  without replacement $\lceil N_S/10 \rfloor$ (where $\lceil x\rfloor$ represents the nearest integer form $x$) states from $S$. Then, for those states, the reward $R(s)$ is drawn randomly and uniformly in $[0,1]$ and for the other states $R(s)=0$. Thus, we obtain a sparse reward function depending only on the states which is the kind of rewards encountered in practice. As we choose finite MDPs, a canonical choice of features $\phi$ is the tabular basis $\phi:S\rightarrow\mathbb{R}^{N_S}$ where $\phi(s)\in \mathbb{R}^{N_S}$ is a vector which is null excepted in $s$ where it is equal to $1$.

A Garnet is a finite MDP where the dynamics $P$ and the reward $R$ is perfectly known. Thus, an optimal policy (playing the role of the expert) can be easily computed by DP (PI in our case). This policy is a key element to build the expert set $D_E$ that fed RCAL and RLED. In our experiments $D_E$ has the following form:
 \begin{equation}
 D_E=(\omega_j)_{\{1\leq j \leq L_E\}},
 \end{equation}
where $\omega_j=(s_{i,j},a_{i,j})_{\{1\leq i \leq H_E\}}$ is a trajectory obtained by starting from a random state $s_{1,j}$ (chosen uniformly in $S$) and applying the expert policy ($\pi_E$) $H_E$ times such that $a_{i,j}=\pi_E(s_{i,j})$ and $s_{i+1,j}=s'_{s_{i,j},a_{i,j}}$. So, $D_E$ is composed by $L_E$ trajectories of $\pi_E$ of length $H_E$ and we have $L_EH_E=N_E$. In addition the data set $D_{RL}$ has the following form:
\begin{equation}
D_{RL}=(\tau_j)_{\{1\leq j \leq L_{RL}\}},
\end{equation}
where $\tau_j=(s_{i,j},a_{i,j},r_{i,j},s_{i,j}'=s'_{s_{i,j},a_{i,j}})_{\{1\leq i \leq H_{RL}\}}$ is a trajectory
obtained by starting from a random state $s_{1,j}$ (chosen uniformly) and applying the the random policy  ($a_{i,j}$ is chosen uniformly from $A$ ) $H_{RL}$ times such that  $s_{i,j}'=s_{i+1,j}$ and $r_{i,j}=R(s_{i,j})$.
So, $D_{RL}$ is composed by $L_{RL}$ trajectories of $\pi_R$ of length $H_{RL}$ and we have $L_{RL}H_{RL}=N_{RL}$. The set $D_{NE}$ corresponds to the set $D_{RL}$ where the reward $r_{i,j}$ is dropped:
\begin{equation}
D_{NE}=(\tau_j)_{\{1\leq j \leq L_{NE}\}},
\end{equation}
where $\tau_j=(s_{i,j},a_{i,j},s_{i,j}'=s'_{s_{i,j},a_{i,j}})_{\{1\leq i \leq H_{NE}\}}$ is a trajectory
obtained by starting from a random state $s_{1,j}$ (chosen uniformly) and applying the the random policy $H_{RL}$ times such that $s_{i,j}'=s_{i+1,j}$.

Finally, as we want to compare gradient descent to DCA for the minimisation of the criteria $J_{RLED}$ and $J_{RCAL}$, it is important to give the parameters of those methods. First, the two methods start form the same starting point which is specified for each experiment. The updates for the gradient descent have the same form as in~Eq.~\eqref{eq: updategrad}. And the number of updates is $100$. To make DCA comparable to gradient descent, we set the number of intermediary convex problems $K$ to $10$ and the number of updates $N$ for the gradient descent of the intermediary problems to $10$. Thus, we have a total of $KN=100$ updates for DCA. In addition, for each gradient descent (it can be the global gradient descent or the one used in the intermediary problems), we set the coefficients $\forall p\in\mathbb{N}^*,\alpha_{p}=1$. 
\subsection{RCAL experiment}
Our first experiment shows the performance improvement of RCAL when $\gamma=0.9$, $\lambda_{RCAL}=0.1$, $H_E$ is increasing and $(L_E,H_{NE},L_{NE})$ are fixed.
It aims at showing that the algorithms perform better when more expert information is available and a fixed amount of knowledge of the dynamics known through $D_{NE}$ is given. It consists in generating $10$ Garnets with $(N_S=100,N_A=5,N_B=1)$ which gives us the set of Garnet problems $\mathfrak{G}=(G_p)_{p=1}^{10}$. On each problem $p$ of the set $\mathfrak{G}$, we compute an optimal and expert policy, $\pi_{E}^p$. The parameter $L_E$ takes its values in the set $(L_E^k)_{k=1}^{10}=(2,4,6,..,20)$ and  $H_E=5$, $H_{NE}=5$, $L_{NE}=20$. Then, for each set of parameters $(L_E^k,H_E,L_{NE},H_{NE})$ and each $G_p$, we compute $20$ expert policy sets $(D_E^{i,p,k})_{i=1}^{20}$ and $20$ random policy sets $(D_{NE}^{i,p,k})_{i=1}^{20}$ which fed the algorithms RCAL and Classif. Overall, we test RCAL on $2000$ sets of data. The starting point of the algorithm RCAL (gradient descent and DCA) and Classif is the null function.

\begin{figure}[h!]
    \centering
    \subfigure[Performance.]{\label{fig1} \includegraphics[width=.48\textwidth]{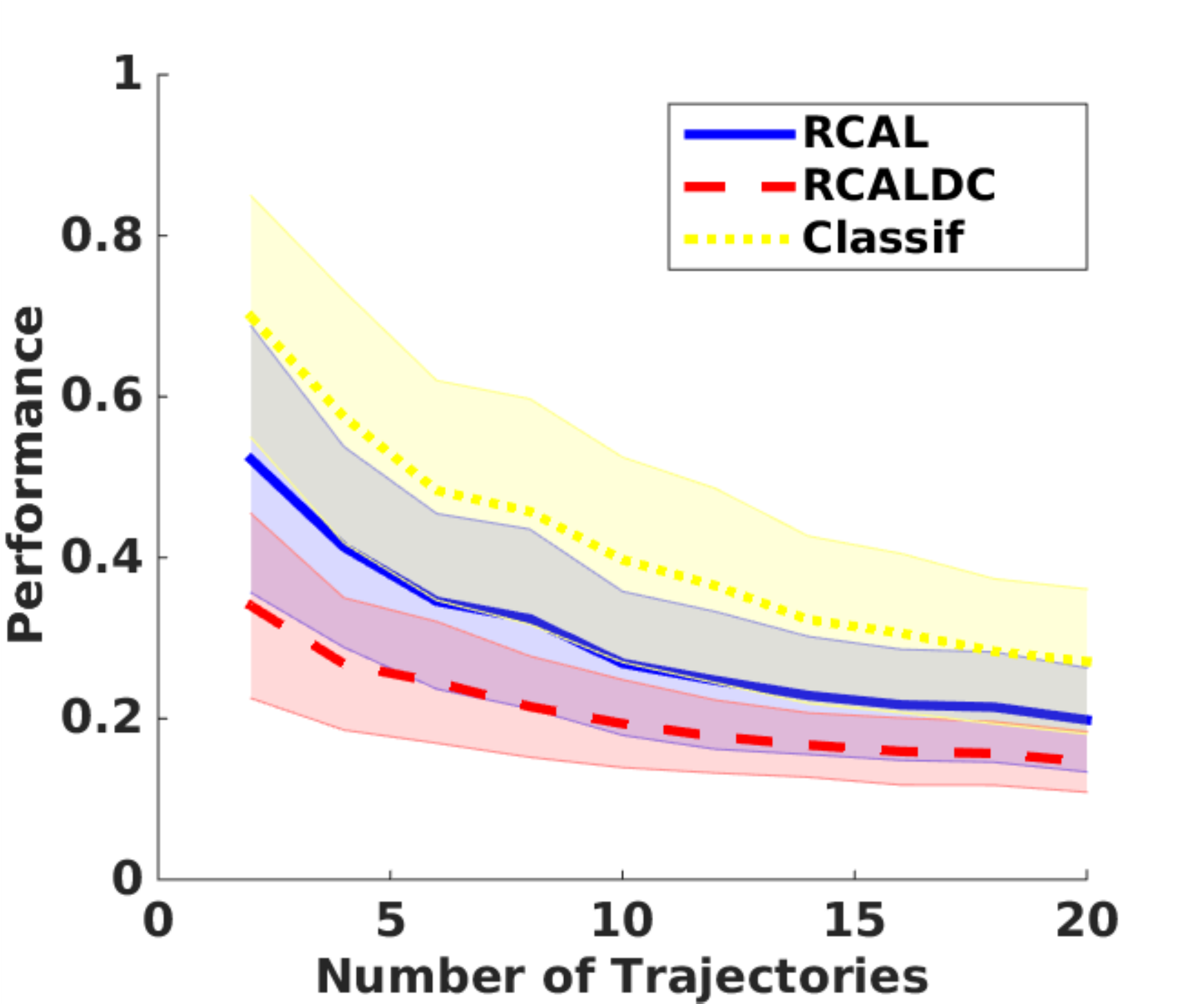}}
    \subfigure[Improvement.]{\label{fig2} \includegraphics[width=.48\textwidth]{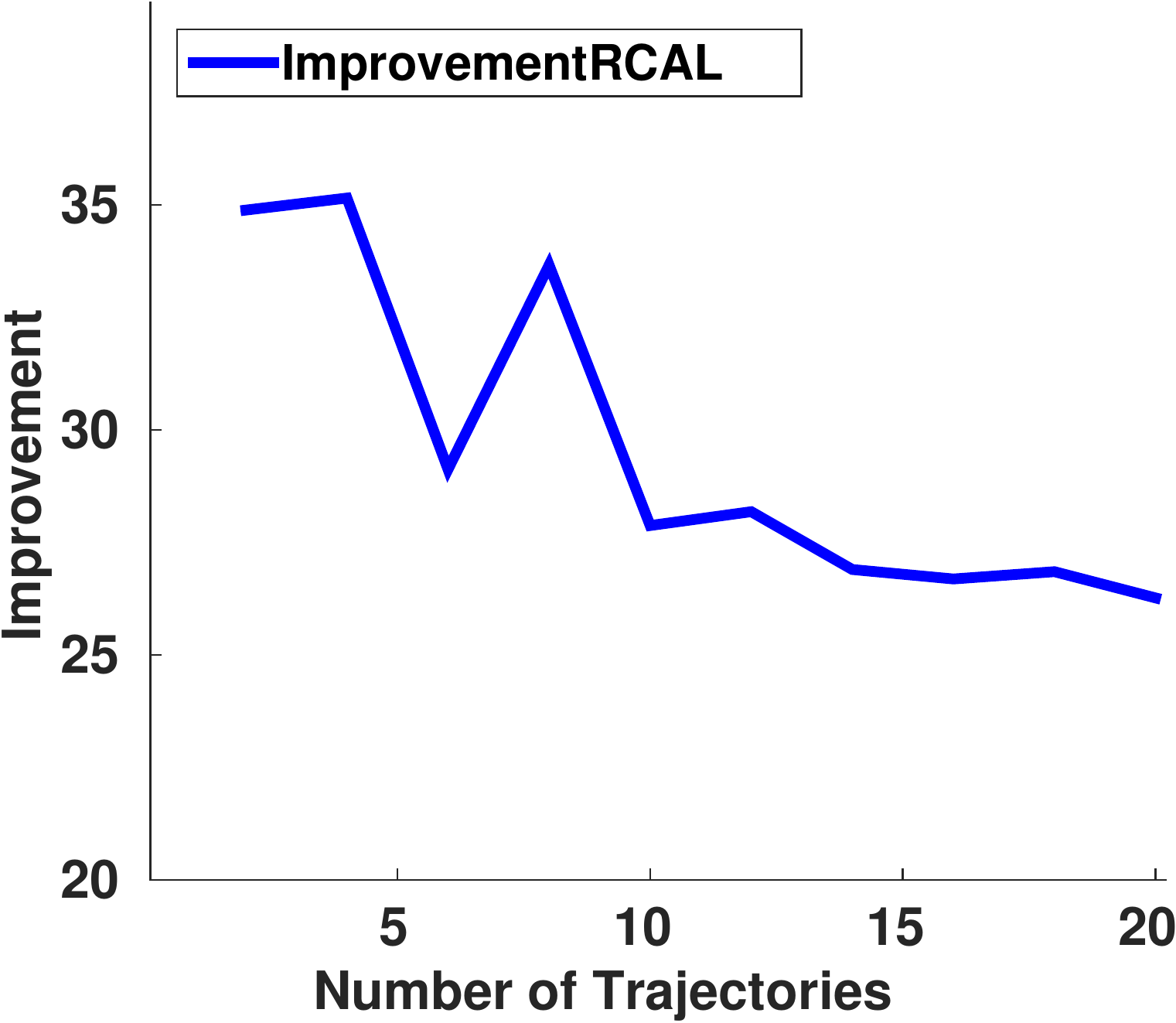}}
    \caption{Garnet Experiment for RCAL.}
    \label{exp1}
\end{figure}

The criterion of performance chosen, for the algorithm $A$ and for each couple $(D_E^{i,p,k},D_{NE}^{i,p,k})$, is the following: 
\begin{equation}
T_{A}^{i,p,k}=\frac{\mathbb{E}_{\rho}[V^{\pi_{E}^p}-V^{\pi_{A}^{i,p,k}}]}{\mathbb{E}_{\rho}[V^{\pi_E^p}]},
\end{equation}
where $\pi_E^p$ is the expert policy, $\pi_{A}^{i,p,k}$ is the policy induced by the algorithm $A$ fed by the couple $(D_E^{i,p,k},D_{NE}^{i,p,k})$ and $\rho$ is the uniform distribution over the state space $S$. For RCAL and Classif, we have $\pi_{A}^{i,p,k}(s)\in\argmax_{a\in A}Q_{\theta^*}(s,a)$ where $\theta^*$ is the output of those algorithms. This criterion of performance is the normalised absolute difference of value-functions between the expert policy and the one induced by the algorithm. Thus, the lesser this criterion is the better. The mean criterion of performance $T^k_A$ for each set of parameters $(L_E^k,H_E,L_{NE},H_{NE})$ is: 
\begin{equation}
T_A^{k}=\frac{1}{200}\sum_{p=1}^{10}\sum_{i=1}^{20} T^{i,p,k}.
\end{equation}
For each algorithm $A$, we plot $(L_E^k,T_A^k)_{k=1}^{10}$ in Fig.~\ref{fig1}. The colored shadows on the figures represent the variance of the algorithms. In order to verify that RCALDC has a better performance than RCAL, we calculate the improvement which is the following ratio:
\begin{equation}
Imp^k=100\frac{T_{RCAL}^{k}-T_{RCALDC}^k}{T_{RCAL}^{k}}.
\end{equation}
This ratio represents in percentage how much RCALDC is better than RCAL. In Fig.~\ref{fig2}, we plot  $(L_E^k,Imp^k)_{k=1}^{10}$.
In Fig.~\ref{fig1} and Fig.~\ref{fig2}, we clearly observe that RCALDC performs in average better than RCAL. In addition, we also compute the number of times that $T_{RCALDC}^{i,p,k}$ is lesser than $T_{RCAL}^{i,p,k}$ over the $2000$ runs of this experiment and we obtain $1732$. Thus, RCALDC is $100\times\frac{1732}{2000}=86.6$\% of the time better than RCAL. Those different elements tend to prove that using DC programming for RCAL improves clearly the performance.

\subsection{RLED experiments}
The second experiment is quite similar to the first except that we use the RLED algorithm. Here $\gamma=0.99$, $\lambda_{RLED}=0.1$, $L_E$ is increasing and $(H_E,H_{RL},L_{RL})$ are fixed. Like the first experiment, it aims at showing that RLED performs better when more expert information is available. It consists in generating a set of Garnets $\mathfrak{G}=(G_p)_{p=1}^{10}$. The parameter $L_E$ takes its values in the set $(L_E^k)_{k=1}^{10}=(1,2,3,..,10)$ and  $H_E=5$, $H_{RL}=5$, $L_{RL}=100$. Then, for each set of parameters $(L_E^k,H_E,L_{RL},H_{RL})$ and each $G_p$, we compute $20$ expert policy sets $(D_E^{i,p,k})_{i=1}^{20}$ and $20$ random policy sets $(D_{RL}^{i,p,k})_{i=1}^{20}$ which fed the algorithms RLED, Classif and LSPI. Here, the starting point of RLED is the ouput of the LSPI algorithm.

\begin{figure}[h!]
    \centering
    \subfigure[Performance.]{\label{fig3} \includegraphics[width=.48\textwidth]{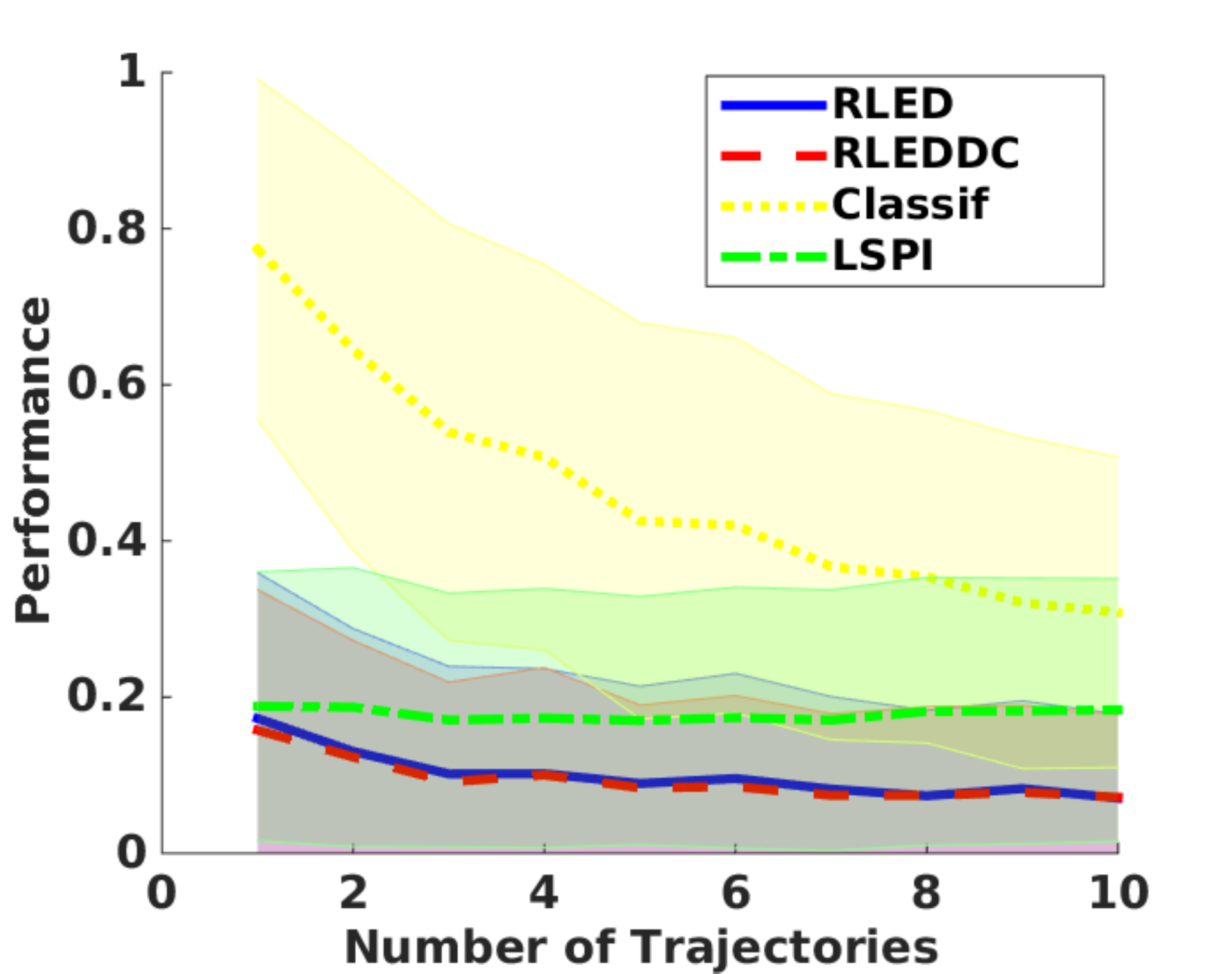}}
    \subfigure[Improvement.]{\label{fig4} \includegraphics[width=.48\textwidth]{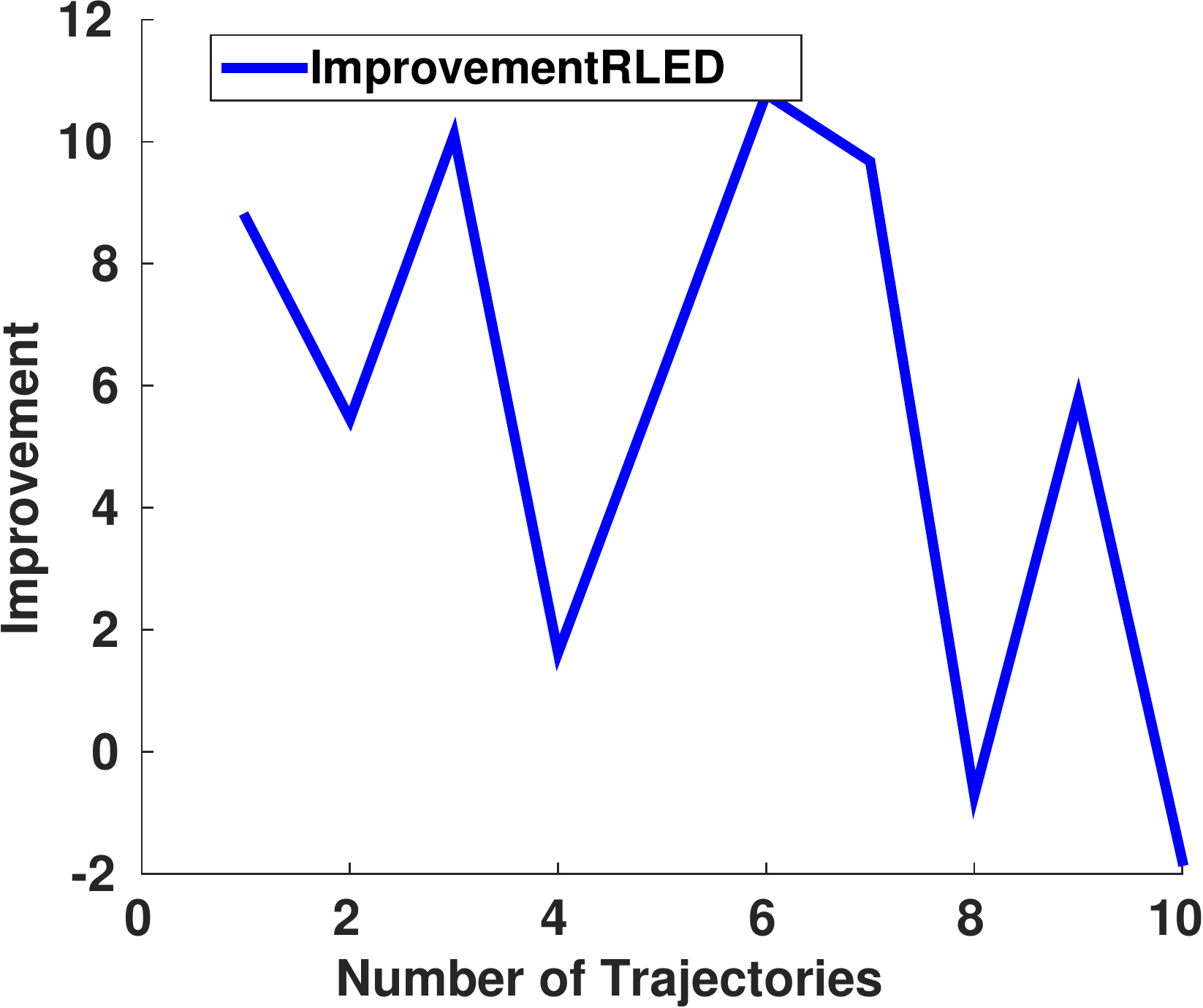}}
    \caption{Garnet Experiment for RLED with $L_E$ increasing.}
    \label{exp2}
\end{figure}

The criterion of performance for the algorithm $A$ and for each couple $(D_E^{i,p,k},D_{RL}^{i,p,k})$ is $T_{A}^{i,p,k}$ which has the same definition as in the first experiment. The mean criterion of performance for each set of parameters $(L_E^k,H_E,L_{RL},H_{RL})$ is $T_A^k$.
For each algorithm $A$, we plot $(L_E^k,T_A^k)_{k=1}^{10}$ in Fig.~\ref{fig3}. In order to verify that RLEDDC has a better performance than RLED, we calculate the improvement which is the following ratio:
\begin{equation}
Imp^k=100\frac{T_{RLED}^{k}-T_{RLEDDC}^k}{T_{RLED}^{k}}.
\end{equation}
This ratio represents in percentage how much RLEDDC is better than RLED. In Fig~\ref{fig4}, we plot  $(L_E^k,Imp^k)_{k=1}^{10}$.
In Fig.~\ref{fig3}, we can not distinguish which algorithm between RLED (with gradient descent) and RLEDDC is better. Thus, we plot a zoom of this curve without the variance to have a better view in Fig.~\ref{fig5}. Even though, it seems that RLEDDC is slightly better, there is no clear difference between the two algorithms. This is principally due to the fact that RLED performs already well without DCA. 

Finally, the third experiment aims at showing that RLED performs better when the information on the model (dynamics and reward) is getting bigger. Here, $\gamma=0.99$ and $\lambda_{RLED}=1$. $\lambda_{RLED}$ is set to a higher value as the set $D_{RL}$ is getting bigger in this experiment and more weight needs to be put on the $J_{RL}$ criterion. The experiment consists in generating a set of Garnets $\mathfrak{G}=(G_p)_{p=1}^{10}$. The parameter $L_{RL}$ takes its values in the set $(L_{RL}^k)_{k=1}^{10}=(50,100,150,..,500)$ and  $L_E=5$, $H_E=5$, $H_{RL}=5$. Then, for each set of parameters $(L_E,H_E,L^k_{RL},H_{RL})$ and each $G_p$, we compute $20$ expert policy sets $(D_E^{i,p,k})_{i=1}^{20}$ and $20$ random policy sets $(D_{RL}^{i,p,k})_{i=1}^{20}$ which fed the algorithms RLED, Classif and LSPI. Here, the starting point of RLED is the ouput of the LSPI algorithm.
Like the previous experiments, for each algorithm $A$, we plot the mean performance $(L_ {RL}^k,T_A^k)_{k=1}^{10}$ in Fig.~\ref{fig5} and  the improvement  $(L_{RL}^k,Imp^k)_{k=1}^{10}$ in Fig.~\ref{fig6}. 
\begin{figure}[h!]
    \centering
    \subfigure[Performance.]{\label{fig5} \includegraphics[width=.48\textwidth]{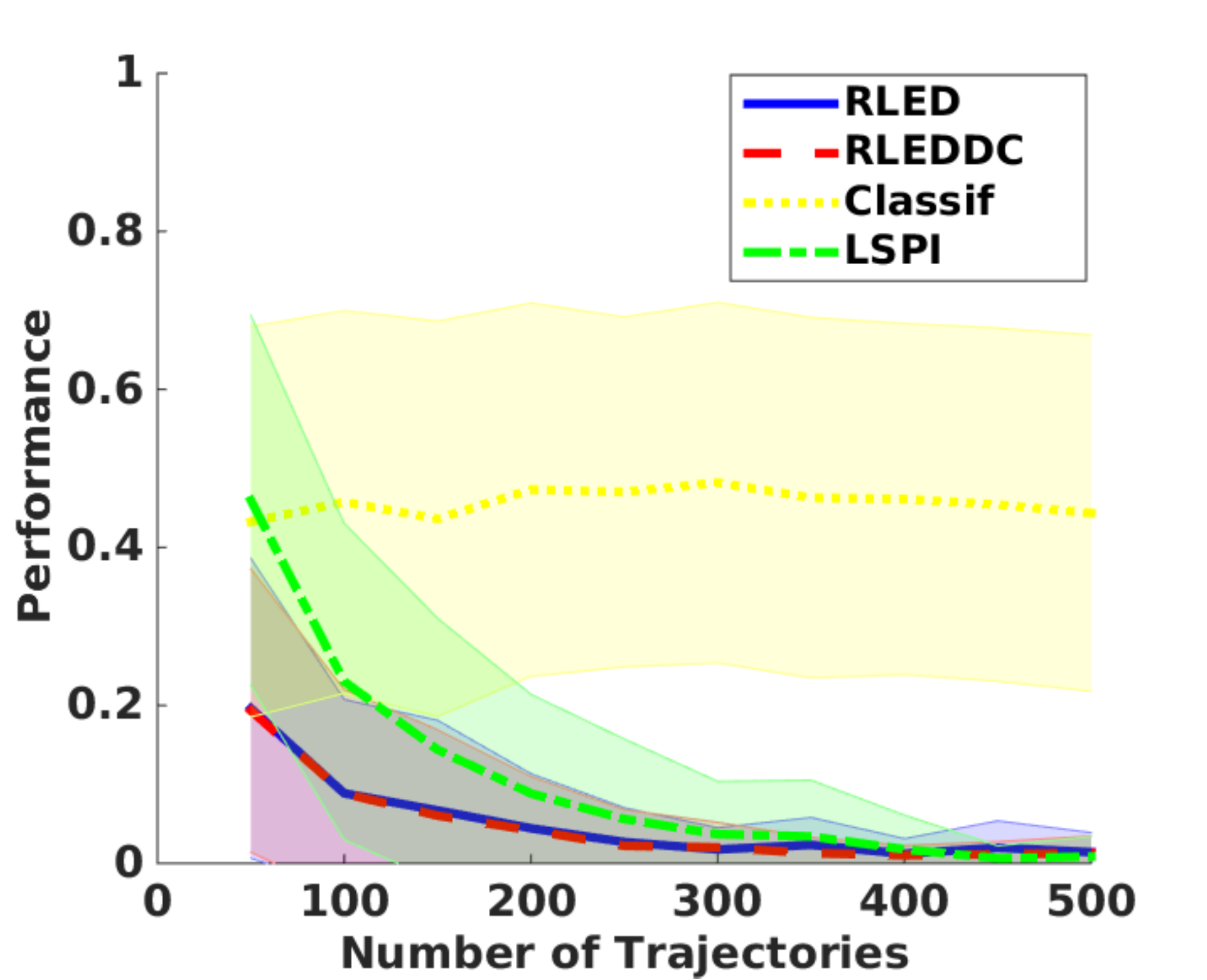}}
    \subfigure[Improvement.]{\label{fig6} \includegraphics[width=.48\textwidth]{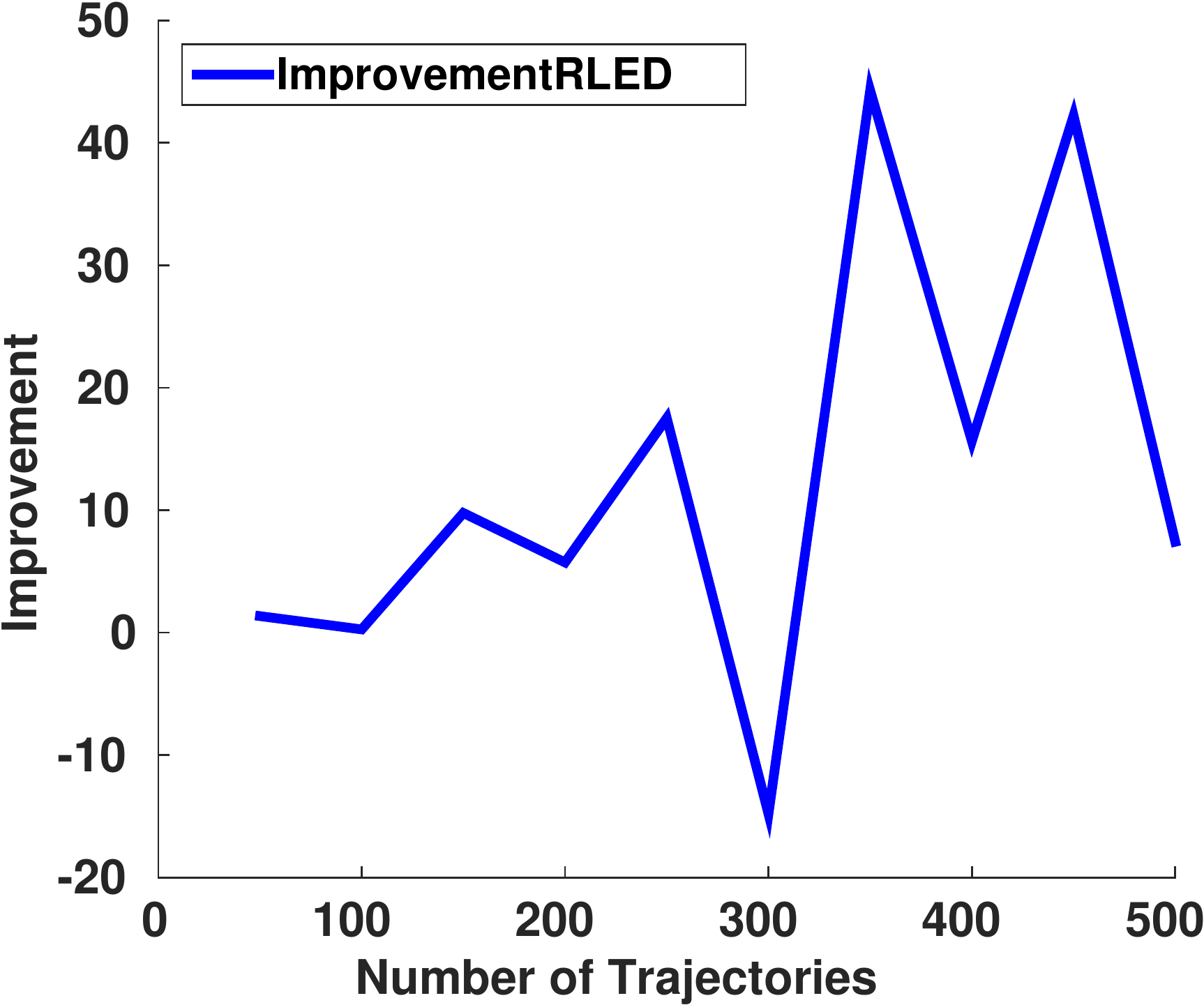}}
    \caption{Garnet Experiment for RLED with $L_{RL}$ increasing.}
    \label{exp3}
\end{figure}
In Fig.~\ref{fig4}, we cannot distinguish between RLED and RLEDDC. A zoom of this plot is proposed in Fig~.\ref{fig6} where there is a slight advantage for RLEDDC. Thus for RLED and contrary to RCAL, even tough there is a slight improvement using DCA, we can not conclude that there is an advantage to use DC programming. 
\begin{figure}[h!]
    \centering
    \subfigure[Zoom of RLED performance with $L_E$ increasing  ]{\label{fig5} \includegraphics[width=.48\textwidth]{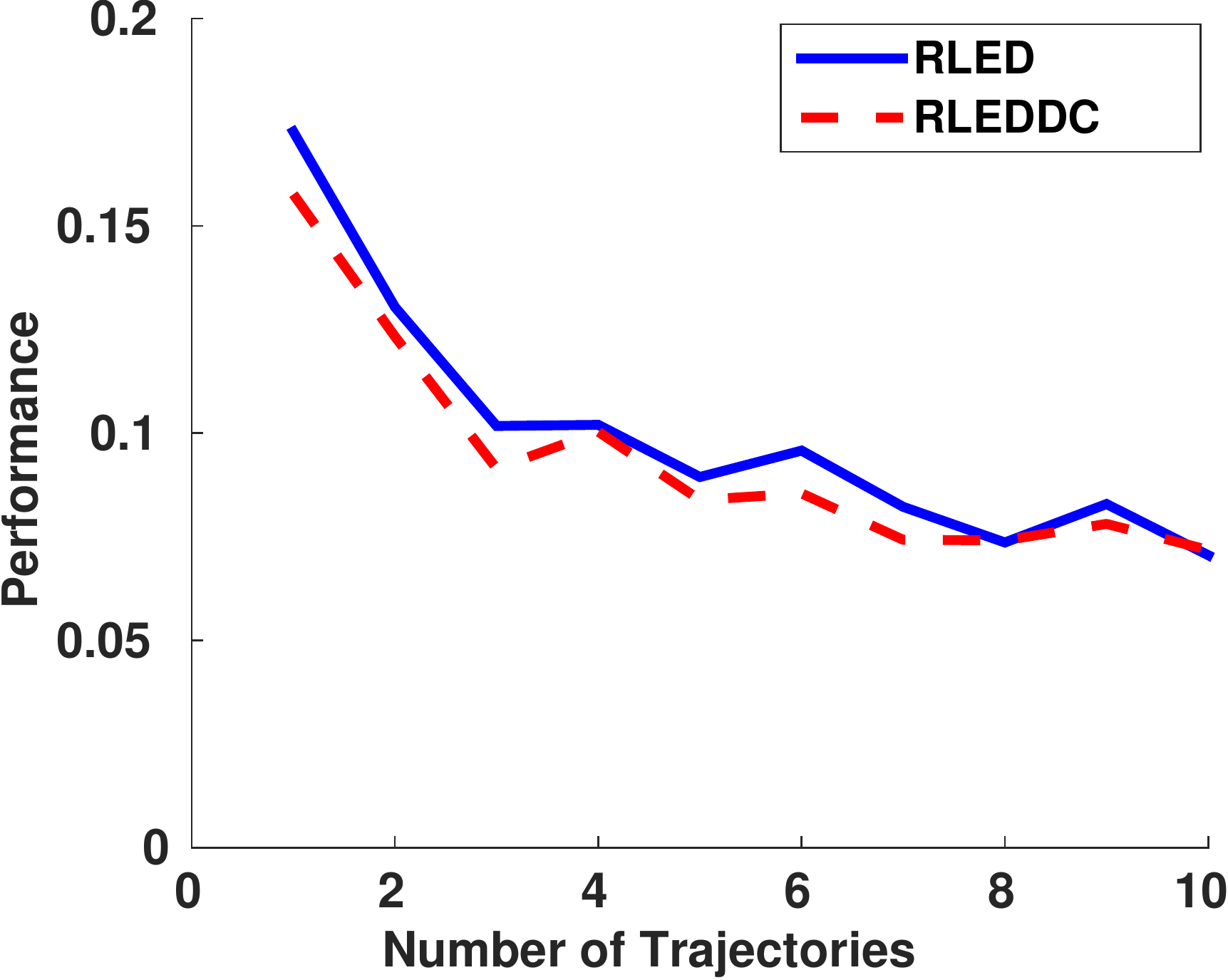}}
    \subfigure[Zoom of RLED performance with $L_{RL}$ increasing]{\label{fig6} \includegraphics[width=.48\textwidth]{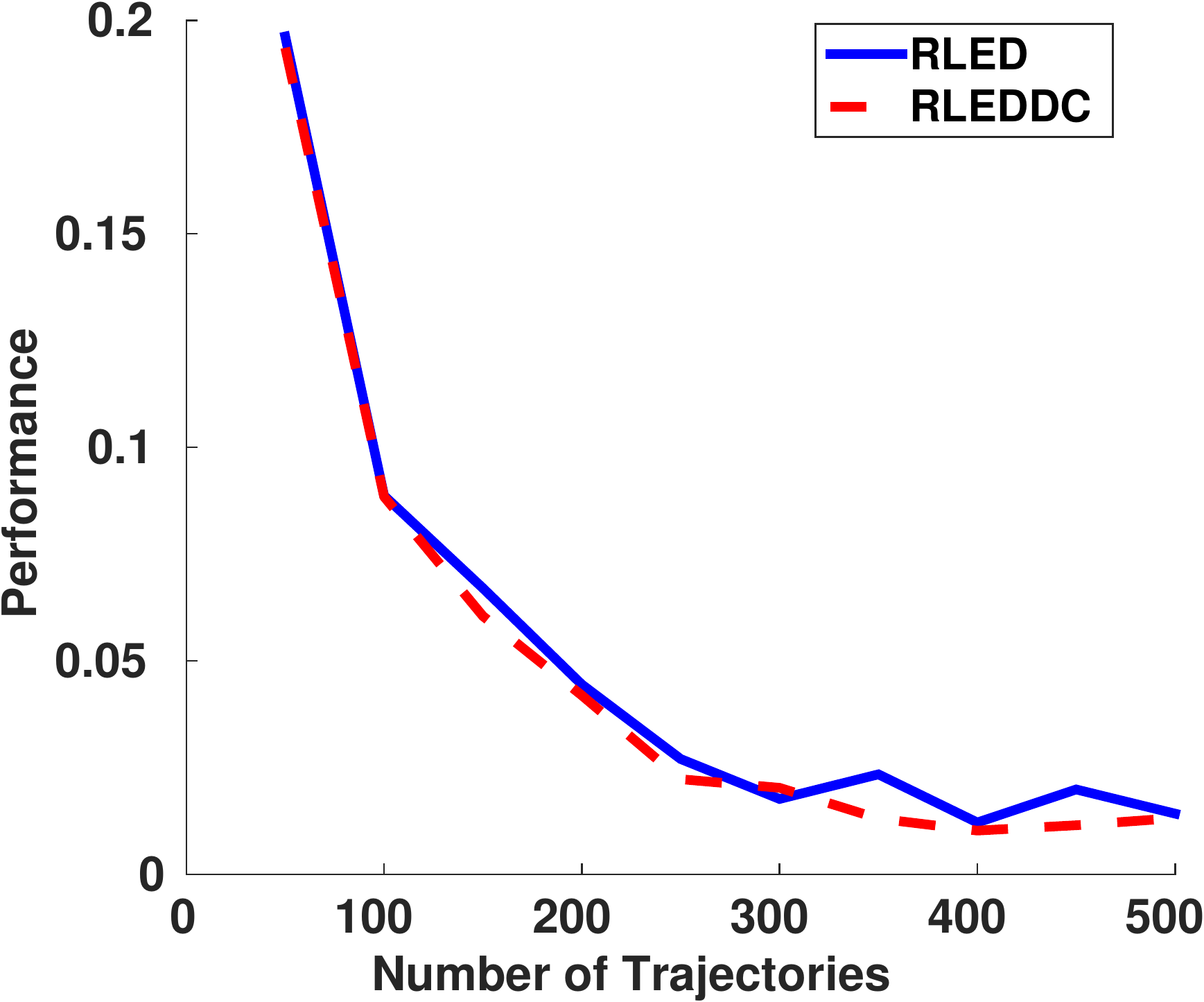}}
    \caption{Zoom for RLED experiments}
    \label{exp3}
\end{figure}
\section{Conclusion and Perspectives} 
\label{section: Conclusion}
In this paper, we showed the implications of seeing the Optimal Bellman Residual (OBR) as a Difference of Convex (DC) functions in the fields of Reinforcement Learning with Expert Demonstrations (RLED) and Learning from Demonstrations (LfD). More precisely, we gave one of the possible DC decompositions of two algorithms, namely RCAL and RLED. In addition, we compared in generic experiments, using randomly constructed Markov Decision Processes (MDPs) called Garnets, the performances of RCAL and RLED using DC programming versus a classical gradient descent. In order to make a fair comparison between the methods, we imposed the same number of updates and the same starting point. Experiments showed a clear advantage of using DC Algorithm (DCA) for the minimisation of the RCAL criterion. However, there was only a slight advantage of DCA for the RLED criterion which can be explained by the fact that RLED with gradient descent already performs well, hence it is quite difficult to improve the method. In conclusion, it seems a promising perspective to use DC programming in fields such as RL and LfD where the goal is to minimise a norm of the OBR which is a DC function. As perspectives, we would like to test several DC decompositions and to start using non-parametric gradient descent in order to solve the intermediary convex problems. Indeed, in RL and more generally in Machine Learning (ML), the choice of features $\phi(s,a)=(\phi_i(s,a))_{i=1}^d$ that represent our hypothesis space is often problem-dependent and need a human expertise. To avoid that step and to automatise even more the algorithm, we want to use non-parametric gradient descent~\citep{grubb2011generalized} which automatically learn its own features. Finally, we would like to use those techniques on large scale and real-life applications to prove their ability to scale-up. 

%\begin{acknowledgements}
%\end{acknowledgements}

\bibliographystyle{spbasic}      % basic style, author-year citations
\bibliography{Biblio_new}   % name your BibTeX data base

\end{document}